\documentclass[11pt]{article}
\usepackage{amssymb}
\usepackage{times}

\title{Hurewicz-like tests for Borel subsets of the plane.\indent}
\author{Dominique LECOMTE}
\date{\it ~Electron. Res. Announc. Amer. Math. Soc.~\rm 11 (2005), 95-102}

\newcommand{\Ana}{{\it\Sigma}^{1}_{1}}
\newcommand{\Ca}{{\it\Pi}^{1}_{1}}
\newcommand{\Borel}{{\it\Delta}^{1}_{1}}

\newcommand{\borel}{{\bf\Delta}^{1}_{1}}
\newcommand{\boraone}{{\bf\Sigma}^{0}_{1}}

\newcommand{\boraxi}{{\bf\Sigma}^{0}_{\xi}}

\newcommand{\bormz}{{\bf\Pi}^{0}_{0}}
\newcommand{\bormone}{{\bf\Pi}^{0}_{1}}
\newcommand{\bormtwo}{{\bf\Pi}^{0}_{2}}

\newcommand{\borapx}{{\bf\Sigma}^{0}_{1+\xi}}

\newcommand{\bormpx}{{\bf\Pi}^{0}_{1+\xi}}
\newcommand{\bormlxi}{{\bf\Pi}^{0}_{<\xi}}

\newcommand{\bormxi}{{\bf\Pi}^{0}_{\xi}}

\newtheorem{thm} {Theorem}
\newtheorem{defi} [thm] {Definition}
\newtheorem{cor} [thm] {Corollary}
\newtheorem{lem} [thm] {Lemma}

\begin{document}

\maketitle

\noindent {\footnotesize {\bf Abstract.} Let $\xi\geq 1$ be a countable ordinal. We study the Borel 
subsets of the plane that can be made $\bormxi$ by refining the Polish topology 
on the real line. These sets are called potentially $\bormxi$. We give a 
Hurewicz-like test to recognize potentially $\bormxi$ sets.}\bigskip\bigskip

\noindent\bf {\Large 1 Preliminaries in dimension one.}\rm\bigskip

 Let us recall some results in dimension one before studying Borel subsets of the 
plane. In descriptive set theory, a standard way to see that a set is complicated 
is to note that it is more complicated than a well-known example. For instance, 
we have the following result (see [SR]):

\begin{thm} (Hurewicz) Let 
$P_{f}\! :=\!\{\alpha\!\in\! 2^{\mathbb{N}} /\exists n\!\in\!
{\mathbb{N}}\ \ \forall m\!\geq\! n\ \ \alpha (m)\! =\! 0\}$, $X$ be a Polish 
space, and $A$ a Borel subset of $X$. Then exactly one of the following holds:\smallskip  

\noindent (a) The set $A$ is $\bormtwo (X)$.\smallskip
  
\noindent (b) There is $u\! :\! 2^{\mathbb{N}}\!\rightarrow\! X$ continuous 
and one-to-one with $P_{f}\! =\! u^{-1}(A)$.\end{thm}

 This result has been generalized to the other Baire classes (see [Lo-SR]). We 
state this generalization in two parts:

\begin{thm} (Louveau-Saint Raymond) Let $\xi\! <\!\aleph_{1}$, 
$A_{1+\xi}\!\in\!\borapx (2^{\mathbb{N}} )$, $X$ be a Polish space, and $A$, $B$ disjoint 
analytic subsets of $X$. One of the following holds:\smallskip
  
\noindent (a) The set $A$ is separable from $B$ by a $\bormpx (X)$ set.\smallskip
  
\noindent (b) There is $u\! :\! 2^{\mathbb{N}}\!\rightarrow\! X$ continuous 
with  $A_{1+\xi}\!\subseteq\! u^{-1}(A)$ and 
$2^{\mathbb{N}}\!\setminus\! A_{1+\xi}\!\subseteq\! u^{-1}(B)$.\smallskip

 If we moreover assume that $A_{1+\xi}\!\notin\!\bormpx$, then this is a dichotomy (in this case, and if $\xi\!\geq\! 2$, then we can have $u$ one-to-one).\end{thm}

\begin{thm} There is a concrete example of 
$A_{1+\xi}\!\in\!\borapx (2^{\mathbb{N}} )\!\setminus\!\bormpx 
(2^{\mathbb{N}} )$, for $\xi\! <\!\aleph_{1}$.\end{thm}

 If we replace $P_{f}$ (resp., $\bormtwo$) with the set $A_{1+\xi}$ given by 
Theorem 3 (resp., $\bormpx$), then we get the generalization of Theorem 1 for 
$\xi\!\geq\! 2$. We state this generalization in two parts for the following 
reasons:\bigskip
 
\noindent $\bullet$ Theorem 2 is valid for any 
$A_{1+\xi}\!\in\!\borapx (2^{\mathbb{N}} )$, and Theorem 1 is of the form 
``There is a typical example such that$\ldots$".\bigskip

\noindent $\bullet$ We will meet again a statement in two parts, in dimension 
two.

\vfill\eject

\noindent\bf {\Large 2 Results with the usual notions of reduction.}\rm\bigskip

 Let us consider the case of dimension two. The usual notion of comparison for 
Borel equivalence relations is the Borel reducibility quasi-order (recall that a 
quasi-order is a reflexive and transitive relation). This means 
that if $X$ (resp., $Y$) is a Polish space, and $E$ (resp., $F$) a Borel 
equivalence relation on $X$ (resp., $Y$), then
$$E\leq_{B}F~\Leftrightarrow ~\exists u\! :\! X\!\rightarrow\! Y\ 
\mbox{\rm Borel\ with}\ E\! =\! (u\!\times\! u)^{-1}(F).$$
Note that this makes sense even if $E$, $F$ are not equivalence 
relations. We will study a natural invariant for $\leq_{B}$. Recall the 
following (see [K]):

\begin{thm} (Kuratowski) Let $X$ be a Polish space, and $(B_{n})$ a 
sequence of Borel subsets of $X$. Then there is a finer Polish topology $\sigma$ 
on $X$ (and thus having the same Borel sets) making the $B_{n}$'s clopen.
\end{thm}

 In particular, if $u\! :\! X\!\rightarrow\! Y$ is Borel, then there is $\sigma$ such 
that $u\! :\! [X,\sigma ]\!\rightarrow\! Y$ is continuous. If moreover 
$E\! =\! (u\!\times\! u)^{-1}(F)$ and $F$ is in some Baire class $\Gamma$, then 
$E\!\in\!\Gamma([X,\sigma ]^2)$. This leads to Definition 5, that  
can be found in [Lo]:

\begin{defi} (Louveau) Let $X$, $Y$ be Polish spaces, $A$ a Borel subset 
of $X\times Y$ and $\Gamma$ a Baire (or Wadge) class. We say that $A$ is 
$potentially~in~\Gamma$ 
${\big(\mbox{denoted}~A\!\in\!\mbox{pot}(\Gamma)\big)}$ iff 
there is a finer Polish topology $\sigma$ (resp., $\tau$) on $X$ (resp., $Y$) such that $A$ is in 
$\Gamma ([X,\sigma ]\!\times\! [Y,\tau ])$.\end{defi}

 The previous result shows that this notion makes sense for product topologies. 
This notion is a natural invariant for $\leq_B$: if $F$ is 
$\mbox{pot}(\Gamma)$ and $E\leq_B F$, then $E$ is $\mbox{pot}(\Gamma)$. 
Using this notion, A. Louveau showed that the collection of $\boraxi$ equivalence 
relations is not cofinal for $\leq_B$, and deduces from this the non existence of 
a maximum Borel equivalence relation for $\leq_B$ (this non existence result is 
due to H. Friedman and L. Stanley). A. Louveau has also more recently noticed that one 
can associate a quasi-order relation $R_{A}\subseteq (X\times 2)^2$ to 
$A\subseteq X^2$ as follows:
$$(x,i)~R_{A}~(y,j)~~\Leftrightarrow ~~(x,i)=(y,j)~\mbox{or}~
[(x,y)\in A~\mbox{and}~(i,j)=(0,1)].$$ 
Using this, one can see that, from the point of view of Borel reducibility, the 
study of Borel quasi-orders is essentially the study of arbitrary Borel subsets 
of the plane. This strengthens the motivation for studying arbitrary Borel 
subsets of the plane, from the point of view of potential complexity. We have a 
result concerning equivalence relations (see [H-K-Lo]):

\begin{thm} (Harrington-Kechris-Louveau) Let $X$ be a Polish space, $E$ a 
Borel equivalence relation on $X$, and $E_{0}\! :=\!\{(\alpha ,\beta )
\!\in\! 2^{\mathbb{N}}\!\times\! 2^{\mathbb{N}} /\exists n\!\in\!{\mathbb{N}}\ \ 
\forall m\!\geq\! n\ \ \alpha (m)\! =\!\beta (m)\}$. Then exactly one of the 
following holds:\smallskip 
 
\noindent (a) The relation $E$ is $\mbox{pot}(\bormone)$.\smallskip  

\noindent (b) $E_{0}\leq_B E$ (with $u$ continuous and one-to-one).\end{thm}

 We will study other structures than equivalence relations (for example 
quasi-orders), and even arbitrary Borel subsets of the plane. We need some other 
notions of comparison. Let $X$, $Y$, $X'$, $Y'$ be Polish spaces, and $A$ (resp., 
$A'$) a Borel subset of $X\!\times\! Y$ (resp., $X'\!\times\! Y'$). We set
$$A\leq^r_{B}A'~\Leftrightarrow ~\exists u\! :\! X\!\rightarrow\! X'\ 
\ \exists v\! :\! Y\!\rightarrow\! Y'\ \ \mbox{Borel\ with}\ A\! =\! 
(u\!\times\! v)^{-1}(A').$$

 We want to extend the previous result to arbitrary Borel subsets of the plane. 
This works partially (see [L1]):

\begin{thm} Let $\Delta (2^{\mathbb{N}} )\! :=\!
\{(\alpha ,\beta )\!\in\! 
2^{\mathbb{N}}\!\times\! 2^{\mathbb{N}}\! /\alpha\! =\!\beta\}$,  
${L_0\! :=\!\{(\alpha ,\beta )\!\in\! 2^{\mathbb{N}}\!\times\! 
2^{\mathbb{N}}\! 
/\alpha\! <_{\mbox{lex}}\!\beta\}}$, $X$, $Y$ be Polish spaces, and $A$ a 
$\mbox{pot}\big(\check D_2(\boraone)\big)$ subset of $X\!\times\! Y$. 
Then exactly one of the following holds:\smallskip 
 
\noindent (a) The set $A$ is $\mbox{pot}(\bormone)$.\smallskip  

\noindent (b) $\neg\Delta (2^{\mathbb{N}} )\leq^r_{B} A$ or 
$L_{0}\leq^r_{B} A$ (with $u$, $v$ continuous and one-to-one).\end{thm}

 Things become more complicated at the level $D_{2}(\boraone)$ 
(differences of two open sets; $\check D_2(\boraone)$ is the dual 
Wadge class of unions of a closed set and of an open set). 

\begin{thm} (a) There is a perfect $\leq^r_{B}$-antichain 
$(A_{\alpha})_{\alpha\in 2^{\mathbb{N}}}\!\subseteq\! 
D_{2}(\boraone )(2^{\mathbb{N}}\!\times\! 2^{\mathbb{N}})$ 
such that $A_{\alpha}$ is $\leq^r_{B}$-minimal among 
$\borel\!\setminus\!\mbox{pot}(\bormone )$ sets, for any 
$\alpha\!\in\! 2^{\mathbb{N}}$.\smallskip 
 
\noindent (b) There is a perfect $\leq_{B}$-antichain 
$(R_{\alpha})_{\alpha\in 2^{\mathbb{N}}}$ such that $R_{\alpha}$ is $\leq_{B}$-minimal 
among $\borel\!\setminus\! \mbox{pot}(\bormone)$ sets, for any 
$\alpha\!\in\! 2^{\mathbb{N}}$. Moreover, $(R_{\alpha})_{\alpha\in 2^{\mathbb{N}}}$ 
can be taken to be a subclass of any of the following classes:\smallskip  

- Graphs (i.e., irreflexive and symmetric relations).\smallskip  

- Oriented graphs  (i.e., irreflexive and antisymmetric relations).\smallskip  

- Quasi-orders.\smallskip    

- Partial orders (i.e., reflexive, antisymmetric and transitive relations).\end{thm}

 In other words, the case of equivalence relations, for which we have a unique 
(up to bi-reducibili-ty) minimal non potentially closed element with Theorem 6, is 
very specific. Theorem 8.(b) says, among other things, that the mixture between 
symmetry and transitivity is very strong.\bigskip 

\noindent\bf Example.\rm ~Let us specify the construction of the antichain in 
(a). We set, for $C\!\subseteq\! 2^{<\mathbb{N}}$,
$$A^{C}\! :=\!\{ (s0\gamma ,s1\gamma )/s\!\in\! C,\gamma\!\in\! 
2^{\mathbb{N}}\}.$$
If $0\!\in\! S\!\subseteq\!\mathbb{N}$ is infinite, then set 
$C_{S}\! :=\!\{ t\!\in\! 2^{<\mathbb{N}}/\hbox{\rm Card}(t)\!\in\! S\}$ 
(where $\hbox{\rm Card}(t)$ is the number of ones in $t$). Such an $S$ is of the 
form $S_{\beta}\! :=\!\{\Sigma_{i<j}~\big(1\! +\!\beta (j)\big)/j\!\in\!
\mathbb{N}\}$, where $\beta\!\in\!\mathbb{N}^{\mathbb{N}}$. 

\begin{thm} The set $A^{C_{S}}$ is minimal for 
${[\borel\!\setminus\!\mbox{pot}(\bormone ),\leq^r_{B}]}$ if
$$\forall p\!\in\!\mathbb{N}\ \ \exists k\!\in\!\mathbb{N}\ \ 
\forall q\!\in\!\mathbb{N}\ \ \exists c\!\in\!\mathbb{N}\!\cap\! 
[q,q\! +\! k]\ \ c\! +\!\big(S\cap [0,p]\big)\! =\! S\cap\big( c\! +\! [0,p]
\big).$$\end{thm}

 It remains to define $\beta_{\alpha}\!\in\! 2^{\mathbb{N}}\!\subseteq\!
\mathbb{N}^{\mathbb{N}}$, for $\alpha\!\in\! 2^{\mathbb{N}}$. We 
inductively define a sequence 
$(s_{\alpha ,n})_{n}\!\subseteq\! 2^{<\mathbb{N}}$ as follows: 
$s_{\alpha ,0}\! :=\! 0$, $s_{\alpha ,1}\! :=\! 1$, 
$s_{\alpha ,n+2}\! :=\! s_{\alpha ,n}^{\alpha (n)+1}s_{\alpha ,n+1}^{\alpha 
(n+1)+1}$. Note that $s_{\alpha ,n}\!\prec_{\not=}\! s_{\alpha ,n+2}$, so that we 
can define $\beta_{\alpha}\! :=\!\hbox{\rm lim}_{n\rightarrow\infty}~
s_{\alpha ,2n}\!\in\! 2^{\mathbb{N}}$. It is suitable: 
$(A^{C_{S_{\beta_{\alpha}}}})_{\alpha\in 2^\mathbb{N}}$ is a perfect 
antichain made of minimal sets for 
${[\borel\!\setminus\!\mbox{pot}(\bormone ),\leq^r_{B}]}$.\smallskip\bigskip

\noindent\bf {\Large 3 Reduction by homomorphism.}\rm\bigskip 

 Theorem 8.(a) shows that the classical notions of reduction (on the whole 
product) don't work, at least at the first level. So we must find another notion 
of comparison.\bigskip

 We have a positive result with another notion, which is in some 
sense ``half of the Borel reducibility ordering". Let $A$ (resp., $A'$) be an 
analytic subset of $X\!\times\! X$ (resp., $X'\!\times\! X'$). We set
$$(X,A)\preceq_{B}(X',A')~\Leftrightarrow ~\exists u\! :\! X\!
\rightarrow\! X'\ \ \mbox{Borel\ with}\ A\!\subseteq\! 
(u\!\times\! u)^{-1}(A').$$
This notion essentially makes sense for irreflexive relations (we can take $u$ to 
be constant if $A'$ is not irreflexive).\bigskip

\noindent\bf Notation.\rm ~Let $\psi\! :\! {\mathbb{N}}\!\rightarrow\! 
2^{<{\mathbb{N}}}$ be the natural bijection (i.e., $\psi (0)\! =\!\emptyset$, 
$\psi (1)\! =\! 0$, $\psi (2)\! =\! 1$, $\psi (3)\! =\! 0^2$, 
$\psi (4)\! =\! 01$, $\psi (5)\! =\! 10$, $\psi (6)\! =\! 1^2$, $\ldots$). Note 
that $|\psi (n)|\!\leq\! n$, so that we can define 
$$s_{n}\! :=\!\psi (n)0^{n-|\psi (n)|}.$$ 
The crucial properties of $(s_{n})$ are that it is dense (there is $n$ such that $t\!\prec\! s_{n}$, for each 
$t\!\in\! 2^{<\omega}$), and that $|s_{n}|\! =\! n$. We put
$$A_{0}\! :=\! A^{\{s_{n}/n\in\mathbb{N}\}}\! =\!
\{ (s_{n}0\gamma ,s_{n}1\gamma )/n\!\in\!{\mathbb{N}},
\gamma\!\in\! 2^{\mathbb{N}}\}.$$
The symmetric set $s(A_0)$ generated by $A_0$ is considered in [K-S-T], where the 
following is essentially proved:

\begin{thm} (Kechris, Solecki, Todor\v cevi\'c) Let $X$ be a Polish 
space, and $A$ an analytic subset of $X\!\times\! X$. Then exactly one of the following 
holds:\smallskip  

\noindent (a) $(X,A)\preceq_{B}({\mathbb{N}},\not=)$.\smallskip  

\noindent (b) $(2^{\mathbb{N}},A_{0})\preceq_{B}(X,A)$ (with $u$ continuous).\end{thm} 

 In [K-S-T], it is conjectured that we can have $u$ one-to-one in Theorem 10.(b). 
This is not the case.\smallskip\bigskip 

\noindent\bf {\Large 4 Reduction on a closed set.}\bigskip\rm

 As a consequence of Theorem 10, we have the following:

\begin{thm} Let $X$, $Y$ be Polish spaces, and $A$ a Borel subset of 
$X\!\times\! Y$. Then exactly one of the following holds:\smallskip  

\noindent (a) The set $A$ is $\mbox{pot}(\bormone)$.\smallskip  

\noindent (b) There are $u\! :\! 2^{\mathbb{N}}\!\rightarrow\! X$ and 
$v\! :\! 2^{\mathbb{N}}\!\rightarrow\! Y$ continuous with 
$A_{0}\! =\! (u\!\times\! v)^{-1}(A)\cap\overline{A_{0}}$.\smallskip  

 Moreover, we can neither ensure that $u$ and $v$ are one-to-one, nor remove 
$\overline{A_{0}}$.\end{thm}  

 So we get a minimum non-potentially closed set if we do not ask for 
a reduction on the whole product. To generalize Theorem 11, the right way to see 
$A_{0}$ seems to be the following. Let $T_{0}$ be the tree associated with 
$\overline{A_{0}}\! =\! A_{0}\cup\Delta (2^{\mathbb{N}})$: 
$$T_{0}\! =\!\{(s,t)\!\in\! 2^{<\mathbb{N}}\!\times\! 
2^{<\mathbb{N}}/s\! =\! t~~\mbox{or}~~\exists n\!\in\!\mathbb{N}~
\exists w\in 2^{<\mathbb{N}}~(s,t)\! =\! (s_{n}0w,s_{n}1w)\}.$$

 The map ${\Delta\! :\! 2^{\mathbb{N}}\!\times\! 2^{\mathbb{N}}\!\rightarrow\! 2^{\mathbb{N}}}$ is the symmetric difference: 
${\Delta (\alpha ,\beta )(i)\! :=\! (\alpha\Delta\beta )(i)\! =\! 1}$ exactly 
when $\alpha (i)\!\not=\!\beta (i)$, for $i\!\in\!\mathbb{N}$. Let 
$S_{1}\! :=\!\{\gamma\!\in\! 2^{\mathbb{N}}/\exists i\!\in\!\mathbb{N}~
\gamma (i)\! =\! 1\}$ be the typical one-dimensional 
$\boraone\!\setminus\!\bormone$ set. We have 
$$A_{0}\! =\!\{ (\alpha ,\beta )\!\in\! 
2^{\mathbb{N}}\!\times\! 2^{\mathbb{N}}/(\alpha ,\beta )\!\in\! [T_{0}]~
\hbox{\rm and}~\alpha\Delta\beta\!\in\! S_{1}\}.$$
This scheme can be generalized. Theorem 11 shows that we cannot have only one 
minimal non potentially closed set, for the reduction on the whole product. The 
reduction is possible on a closed set (the closure of $A_{0}$). This closure does 
not explain why we cannot have a reduction on the whole product. This comes from 
Theorem 8. The orthogonality between the examples appearing in the antichains of 
its statement comes from different types of cycles. This will give a better 
explanation than the closure. We will replace the closure with a closed set, that 
will be seen as the set of branches of some tree on $2\!\times\! 2$. This tree 
will have the acyclicity properties that we need. This leads to the following 
definition:

\begin{defi} Let $R$ be a relation on a set $E$.\smallskip 

\noindent $\bullet$ An $R\mbox{-}path$ is a finite sequence 
$(e_{i})_{i\leq n}\!\subseteq\! E$ such that $(e_{i},e_{i+1})\!\in\! R$, for 
$i\! <\! n$.\smallskip

\noindent $\bullet$ An $R\mbox{-}cycle$ is an $R$-path 
$(e_{i})_{i\leq n}$ such that $n\!\geq\! 3$ and
$$[0\!\leq\! i\!\not=\! j\!\leq\! n~\mbox{and}~e_{i}\! =\! 
e_{j}]~\Leftrightarrow ~\{ i,j\}\! =\!\{ 0,n\}.$$
$\bullet$ We say that $R$ is $acyclic$ if there is no 
$R$-cycle.\smallskip

\noindent $\bullet$ We say that a tree $T$ on $2\!\times\! 2$ is 
$uniformly~acyclic$ if, for each $p\! >\! 0$,\smallskip  

\noindent (a) The relation $T\cap (2^p\!\times\! 2^p)$ is irreflexive and 
antisymmetric.\smallskip  

\noindent (b) The symmetric relation $s\big( T\cap (2^p\!\times\! 2^p)\big)$ 
generated by $T\cap (2^p\!\times\! 2^p)$ is acyclic.\end{defi}

 The main new results in this paper are the following:

\begin{thm} (Debs-Lecomte) Let $T$ be a uniformly acyclic tree, 
$\xi\! <\!\aleph_{1}$, $A_{1+\xi}$ in $\borapx ([T])$, $X$, $Y$ Polish 
spaces, and $A$, $B$ disjoint analytic subsets of $X\!\times\! Y$. Then one of 
the following holds:\smallskip  

\noindent (a) The set $A$ is separable from $B$ by a $\mbox{pot}(\bormpx )$ set.\smallskip  

\noindent (b) There are $u\! :\! 2^{\mathbb{N}}\!\rightarrow\! X$ and 
$v\! :\! 2^{\mathbb{N}}\!\rightarrow\! Y$ continuous such that the inclusions  
$A_{1+\xi}\!\subseteq\! (u\!\times\! v)^{-1}(A)$ and 
${[T]\!\setminus\! A_{1+\xi}\subseteq (u\times v)^{-1}(B)}$ hold.\smallskip

 If we moreover assume that $A_{1+\xi}\!\notin\!\mbox{pot}(\bormpx )$, then this is a dichotomy.\end{thm}

 This result has initially been shown by D. Lecomte when $1\! +\!\xi$ is a 
successor ordinal. Then G. Debs proved it when $1\! +\!\xi$ is a limit ordinal. 
The proof of Theorem 13 uses the representation Theorem for Borel sets in [D-SR]. 
Notice that we can deduce Theorem 2 from the proof of Theorem 13. Theorem 13 is 
the analog of Theorem 2 in dimension two (see [Lo-SR], also Theorem III-2.1 in 
[D-SR]). The tree $T$ has to be small enough, since there is no possibility to 
have a reduction on the whole product. But as the same time, $T$ has to be big 
enough to ensure the existence of complicated sets inside $[T]$:

\begin{thm} There are concrete examples of:\smallskip

\noindent (a) A uniformly acyclic tree $T$.\smallskip

\noindent (b) A set $A_{1+\xi}\!\in\!\borapx ([T])\!\setminus\!\mbox{pot}(\bormpx )$, 
for $\xi\! <\!\aleph_{1}$.\end{thm}

 This result is the complement of Theorem 13 (which is true with 
$T\! :=\!\emptyset$!). Again, the couple Theorems 13-14 is the analog of 
the couple Theorems 2-3.

\vfill\eject

\noindent\bf {\Large 5 The examples.}\rm\bigskip

 Let us specify the examples of Theorem 14. Let 
$\varphi\! =\! (\varphi_{0},\varphi_{1})\!:\! {\mathbb{N}}\!\rightarrow
\! {\mathbb{N}}^2$ be the natural bijection. More precisely, we set, for $q\!\in\!\mathbb{N}$,
$$M(q)\! :=\!\mbox{max}
\{ m\!\in\!\mathbb{N}/\Sigma_{k\leq m}~k\!\leq q\}.$$
Then we define $\varphi (q)\! =\!\big(\varphi_{0}(q),\varphi_{1}(q)\big)\! :=\!
\big(M(q)\! -\! q\! +\! (\Sigma_{k\leq M(q)}~k),
q\! -\! (\Sigma_{k\leq M(q)}~k)\big)$. One can check that 
$<i,j>:=\!\varphi^{-1}(i,j)\! =\! (\Sigma_{k\leq i+j}~k)\! +\! j$. More 
concretely, we get 
$$\varphi [\mathbb{N} ]\! =\!
\{(0,0);(1,0);(0,1);(2,0);(1,1);(0,2);\ldots\}.$$

\begin{defi} We say that 
$E\!\subseteq\!\bigcup_{q\in{\mathbb{N}}}~2^q\!\times\! 2^q$ is a $test$ if\smallskip  

\noindent (a) $\forall q\!\in\!{\mathbb{N}}\ \ \exists ! (s_{q},t_{q})\!\in\! E
\cap (2^q\!\times\! 2^q)$.\smallskip  
 
\noindent (b) $\forall m,q\!\in\!{\mathbb{N}}\ \ \forall u\!\in\! 
2^{<{\mathbb{N}}}\ \ \exists v\!\in\! 2^{<{\mathbb{N}}}\ 
\ (s_{q}0uv,t_{q}1uv)\!\in\! E$ and $\varphi_{0}(|t^{}_{q}1uv|\! -\! 1)\! =\! m$.
\smallskip  

\noindent (c) $\forall n\! >\! 0\ \ \exists q\! <\! n\ \ \exists w\!\in\! 
2^{<{\mathbb{N}}}\ \ s_{n}\! =\! s_{q}0w$ and $t_{n}\! =\! t_{q}1w$.\smallskip  

 We will call $T$ the tree generated by a test 
$E\! =\!\{(s_{q},t_{q})/q\!\in\!{\mathbb{N}}\}$:
$$T\! :=\!\{(s,t)\!\in\! 2^{<{\mathbb{N}}}\!\times\! 2^{<{\mathbb{N}}}/
s\! =\! t\! =\!\emptyset\ ~\mbox{or}~\ \exists q\!\in\!{\mathbb{N}}\ \ 
\exists w\!\in\! 2^{<{\mathbb{N}}}\ \ s\! =\! s_{q}0w\ ~\mbox{and}~\ 
t\! =\! t_{q}1w\}.$$\end{defi}

 One can show the existence of a test, and that $T$ is uniformly acyclic. The 
uniqueness condition in (a) and condition (c) ensure that $T$ is small enough, 
and also the acyclicity. The existence condition in (a) and condition (b) ensure 
that $T$ is big enough. More specifically, if $X$ is a Polish space and $\sigma$ 
a finer Polish topology on $X$, then there is a dense $G_{\delta}$ subset of $X$ 
on which the two topologies coincide. The first part of condition (b) ensures the 
possibility to get inside the square of a dense $G_{\delta}$ subset of $2^\omega$. 
The examples of Theorem 14.(b) are constructed using the examples in [Lo-SR]. 
Conditions on the verticals appear, and the second part of condition (b) gives a 
control on the choice of verticals.\bigskip

\noindent\bf Notation.\rm ~In [Lo-SR], Lemma 3.3, the map 
$\rho_{0}\! :\! 2^{\mathbb{N}}\!\rightarrow\! 2^{\mathbb{N}}$ defined as 
follows is introduced:
$$\rho_{0}(\varepsilon )(i)\! :=\!\left\{\!\!\!\!\!\!
\begin{array}{ll}
& 1\ \mbox{if}\ \varepsilon\big(\! <i,j>\!\big)\! =\! 0,\ 
\hbox{\rm for\ each\ }j\!\in\!\mathbb{N}\mbox{,}\cr 
& 0\ \mbox{otherwise.}
\end{array}
\right.$$
In this paper, 
$\rho_{0}^{\xi}\! :\! 2^{\mathbb{N}}\!\rightarrow\! 2^{\mathbb{N}}$ is 
also defined for $\xi\! <\!\aleph_{1}$ as follows, by induction on $\xi$ (see the 
proof of Theorem 3.2). We put 
${\rho^0_{0}\! :=\!\hbox{\rm Id}_{2^{\mathbb{N}}}}$, 
$\rho^{\eta +1}_{0}\!:=\!\rho^{}_{0}\circ\rho^{\eta}_{0}$. If $\lambda\! >\! 0$ 
is limit, then fix $(\xi^\lambda_{k})\!\subseteq\!\lambda\!\setminus\!\{ 0\}$ 
such that $$\Sigma_{k}~\xi^\lambda_{k}\! =\!\lambda .$$ 
For $\varepsilon \!\in\! 2^{\mathbb{N}}$ and $k\!\in\!\mathbb{N}$, we define 
$(\varepsilon )^k\!\in\! 2^{\mathbb{N}}$ 
by $(\varepsilon )^k(i)\! :=\!\varepsilon (i\! +\! k)$. We also define 
$\rho^{(k,k+1)}_{0}\! :\! 2^{\mathbb{N}}\!\rightarrow\! 2^{\mathbb{N}}$ by 
$$\rho^{(k,k+1)}_{0}(\varepsilon )(i)\! :=\!\left\{\!\!\!\!\!\!
\begin{array}{ll} 
& \varepsilon (i)~\mbox{if}~i\! <\! k\mbox{,}\cr 
& \rho_{0}^{\xi^\lambda_{k}}\big( (\varepsilon )^k\big)(i\! -\! k)~
\mbox{if}~i\!\geq\! k.
\end{array}
\right.$$
We set $\rho^{(0,k+1)}_{0}\! :=\!\rho^{(k,k+1)}_{0}\circ\rho^{(k-1,k)}_{0}\circ
\ldots\circ\rho^{(0,1)}_{0}$ and 
$\rho^\lambda_{0}(\varepsilon )(k)\! :=\!\rho^{(0,k+1)}_{0}(\varepsilon )(k)$.\bigskip

 The set ${H_{1+\xi}\!:=\! (\rho^{\xi}_{0})^{-1}(\{ 0^\infty\})}$ is also 
introduced, and the authors show that $H_{1+\xi}$ is $\bormpx\setminus\borapx$ 
(see Theorem 3.2).

\vfill\eject

\noindent $\bullet$ The map $S:2^{\mathbb{N}}\rightarrow 2^{\mathbb{N}}$ 
is the shift map: $S(\alpha )(i)\! :=\!\alpha (i\! +\! 1)$.\bigskip

\noindent $\bullet$ Let $T$ be the tree generated by a test. We put, 
for $\xi\! <\!\aleph_{1}$,
$$A_{1+\xi}\! :=\!\{ (\alpha ,\beta )\!\in\! 2^{\mathbb{N}}\!
\times\! 2^{\mathbb{N}}/(\alpha ,\beta )\!\in\! [T]\ \mbox{and}
\ S(\alpha\Delta\beta )\!\notin\! H_{1+\xi}\}.$$
Then $A_{1+\xi}$ is $\borapx ([T])\!\setminus\!\mbox{pot}(\bormpx )$. We 
introduce a notation to state the crucial Lemma used to show it:\bigskip

\noindent\bf Notation.\rm ~We define $p\! :\!\mathbb{N}^{<\mathbb{N}}\!
\setminus\!\{\emptyset\}\!\rightarrow\!\mathbb{N}$. We 
actually define $p(s)$ by induction on $|s|$:
$$p(s)\! :=\!\left\{\!\!\!\!\!\!
\begin{array}{ll}
& s(0)~\mbox{if}~|s|\! =\! 1\mbox{,}\cr 
& <p\big(s\lceil (|s|\! -\! 1)\big),s(|s|\! -\! 1)>~\mbox{otherwise.}
\end{array}
\right.$$
Notice that $p\vert_{\mathbb{N}^n}\! :\!\mathbb{N}^n\!\rightarrow\!
\mathbb{N}$ is a bijection for each $n\!\geq\! 1$.

\begin{lem} Let $G$ be a dense $G_{\delta}$ subset of 
$2^{\mathbb{N}}$. Then there are $\alpha_{0}\!\in\! G$ and 
$f\! :\! 2^{\mathbb{N}}\!\rightarrow\! G$ continuous such that, for each 
$\alpha\!\in\! 2^{\mathbb{N}}$,\smallskip

\noindent (a) $\big(\alpha_{0},f(\alpha )\big)\!\in\! [T]$.\smallskip

\noindent (b) For each $t\!\in\!\mathbb{N}^{<\mathbb{N}}$, and each 
$m\!\in\!\mathbb{N}$,\smallskip 

(i) $\alpha\big( p(tm)\big)\! =\! 1~\Rightarrow ~\exists m'\!\in\!\mathbb{N}
\ \ \big(\alpha_{0}\Delta f(\alpha )\big)\big( p(tm')\! +\! 1\big)\! =\! 1$.
\smallskip 

(ii) $\big(\alpha_{0}\Delta f(\alpha )\big)\big( p(tm)\! +\! 1\big)\! =\! 1~
\Rightarrow ~\exists m'\!\in\!\mathbb{N}\ \ \alpha\big( p(tm')\big)\! =\! 1$.
\end{lem}\smallskip

\noindent\bf {\Large 6 Complements of the main results.}\rm\bigskip

 Now we come to consequences of Theorems 13 and 14. To state them, we need some 
more notation. We use some tools from effective descriptive set theory 
(the reader should see [M] for basic notions about it).\bigskip

\noindent\bf Notation.\rm ~Let $X$ be a recursively presented Polish space. We 
denote by ${\it\Delta}_{X}$ the topology on $X$ generated by $\Borel (X)$. This 
topology is Polish (see the proof of Theorem 3.4 in [Lo]). We set 
$$\tau_{1}\!:=\! {\it\Delta}_{X}\!\times\! {\it\Delta}_{Y}$$ 
if $Y$ is also a recursively presented Polish space.\bigskip

\noindent $\bullet$ Let $2\!\leq\!\xi\! <\!\omega^{\hbox{\rm CK}}_{1}$. The 
topology $\tau_{\xi}$ is generated by 
${\Ana (X\!\times\! Y)\cap\bormlxi (\tau_{1})}$. Note that 
$${\boraone(\tau_{\xi})\!\subseteq\!\boraxi (\tau_{1})}\mbox{,}$$ 
so that ${\bormone(\tau_{\xi})\!\subseteq\!\bormxi (\tau_{1})}$.\bigskip

\noindent $\bullet$ Recall the existence of $\Ca$ sets 
$W^{X}\!\subseteq\!\mathbb{N}$, 
$C^{X}\!\subseteq\!\mathbb{N}\!\times\! X$ with 
${\Borel (X)\! =\!\{ C^{X}_{n}/n\!\in\! W^{X}\}}$, 
$${\{(n,x)\!\in\!\mathbb{N}\!\times\!\! X/n\!\in\! W^{X}\hbox{\rm and}~
x\!\notin\! C_{n}^{X}\}\!\in\!\Ca (\mathbb{N}\!\times\! X)}$$ 
(see [H-K-Lo], Theorem 3.3.1).

\vfill\eject

\noindent $\bullet$ Set 
$\mbox{pot}(\bormz)\! :=\!\borel (X)\!\times\!\borel (Y)$ and, for 
$\xi\! <\!\omega^{\mbox{CK}}_{1}$,
$$W^{X\times Y}_{\xi}:=\{ p\!\in\! W^{X\times Y}/
C^{X\times Y}_{p}\!\in\! \mbox{pot}(\bormxi)\}.$$
We also set ${W^{X\times Y}_{<\xi}\! :=\!\bigcup_{\eta <\xi}~W^{X\times Y}_{\eta}}$.

\begin{thm} (Debs-Lecomte-Louveau) Let $T$ given by Theorem 14, 
$\xi\! <\!\omega^{\mbox{CK}}_{1}$, $A_{1+\xi}$ given by Theorem 14, and $X$, 
$Y$ be recursively presented Polish spaces.\smallskip

\noindent (1) Let $A$, $B$ be disjoint $\Ana$ subsets of $X\!\times\! Y$. The 
following are equivalent:\smallskip

\noindent (a) The set $A$ cannot be separated from $B$ by a $\mbox{pot}(\bormpx )$ set.\smallskip

\noindent (b) The set $A$ cannot be separated from $B$ by a 
$\Borel\cap\mbox{pot}(\bormpx )$ set.\smallskip

\noindent (c) The set $A$ cannot be separated from $B$ by a $\bormpx (\tau_1)$ set.\smallskip

\noindent (d) $\overline{A}^{\tau_{1+\xi}}\cap B\!\not=\!\emptyset$.\smallskip

\noindent (e) There are $u\! :\! 2^{\mathbb{N}}\!\rightarrow\! X$ and 
$v\! :\! 2^{\mathbb{N}}\!\rightarrow\! Y$ continuous such that the inclusions  
$A_{1+\xi}\!\subseteq\! (u\!\times\! v)^{-1}(A)$ and 
$[T]\!\setminus\! A_{1+\xi}\!\subseteq\! (u\!\times\! v)^{-1}(B)$ hold.\smallskip

\noindent (2) The sets $W^{X\times Y}_{1+\xi}$ and $W^{X\times Y}_{<1+\xi}$ are $\Ca$.
\end{thm}

 The equivalence between (a), (b) and (c), and also (2), are proved in [Lo]. We 
can assume this equivalence and (2), then prove Theorems 13, 14, and then prove 
Theorem 17. We can also prove directly Theorem 17 by induction on $\xi$. An 
immediate consequence of this is the following, proved in [Lo]:

\begin{cor} (Louveau) Let $\xi\! <\!\omega^{\hbox{\rm CK}}_{1}$, 
$X$, $Y$ be recursively presented Polish spaces, and $A$ a $\Borel$ subset of 
$X\!\times\! Y$. The following are equivalent:\smallskip

\noindent (a) The set $A$ is $\mbox{pot}(\bormpx )$.\smallskip

\noindent (b) The set $A$ is $\bormpx (\tau_1)$.\end{cor}

 We also have the following consequence of Theorems 13 and 14:

\begin{cor} (Debs-Lecomte) Let $\xi\! <\!\aleph_1$. There is a 
Borel subset $A_{1+\xi}$ of $2^{\mathbb{N}}\!\times\! 2^{\mathbb{N}}$ 
such that for any Polish spaces $X$, $Y$, and for any disjoint 
analytic subsets $A$, $B$ of $X\!\times\! Y$, exactly one of the 
following holds:\smallskip

\noindent (a) The set $A$ is separable from $B$ by a $\mbox{pot}(\bormpx )$ set.\smallskip

\noindent (b) There are $u\! :\! 2^{\mathbb{N}}\!\rightarrow\! X$ and 
$v\! :\! 2^{\mathbb{N}}\!\rightarrow\! Y$ continuous such that the inclusions  
$A_{1+\xi}\!\subseteq\! (u\!\times\! v)^{-1}(A)$ and 
$\overline{A_{1+\xi}}\!\setminus\! A_{1+\xi}\!\subseteq\! (u\!\times\! v)^{-1}(B)$ hold.\smallskip

 Moreover we can neither ensure that $u$ and $v$ are one-to-one if 
$\xi\!\leq\! 1$, nor replace $\overline{A_{1+\xi}}\!\setminus\! A_{1+\xi}$ with 
$(2^{\mathbb{N}}\!\times\! 2^{\mathbb{N}})\!\setminus\! A_{1+\xi}$.
\end{cor}

 The one-to-one complement is due to D. Lecomte (see Theorem 11 when $\xi =0$, 
and Theorem 15 in [L2] when $\xi =1$). The latter complement has initially been 
shown by D. Lecomte when $\xi\leq 1$ (see for example Theorem 11). Then G. Debs 
found a simpler proof, which moreover works in the general case.\bigskip\smallskip 

\noindent\bf {\Large 7 References.}\rm\bigskip

\noindent [D-SR]\ \ G. Debs and J. Saint Raymond,~\it Borel liftings of Borel sets: some decidable and undecidable statements,~\rm Mem. Amer. Math. Soc.~187, 876 (2007)

\noindent [H-K-Lo]\ \ L. A. Harrington, A. S. Kechris and A. Louveau,~\it A Glimm-Effros 
dichotomy for Borel equivalence relations,~\rm J. Amer. Math. Soc.~3 (1990), 903-928

\noindent [K]\ \ A. S. Kechris,~\it Classical Descriptive Set Theory,~\rm 
Springer-Verlag, 1995

\noindent [K-S-T]\ \ A. S. Kechris, S. Solecki and S. Todor\v cevi\'c,~\it Borel chromatic numbers,\ \rm 
Adv. Math.~141 (1999), 1-44

\noindent [L1]\ \ D. Lecomte,~\it Classes de Wadge potentielles et th\'eor\`emes d'uniformisation 
partielle,~\rm Fund. Math.~143 (1993), 231-258

\noindent [L2]\ \ D. Lecomte,~\it Complexit\'e des bor\'eliens~\`a coupes 
d\'enombrables,~\rm Fund. Math.~165 (2000), 139-174

\noindent [Lo]\ \ A. Louveau,~\it Ensembles analytiques et bor\'eliens dans les espaces produit,~\rm 
Ast\'erisque (S. M. F.) 78 (1980)

\noindent [Lo-SR]\ \ A. Louveau and J. Saint Raymond,~\it Borel classes and closed games : Wadge-type and Hurewicz-type results,~\rm Trans. A. M. S.~304 (1987), 431-467

\noindent [M]\ \ Y. N. Moschovakis,~\it Descriptive set theory,~\rm North-Holland, 1980

\noindent [SR]\ \ J. Saint Raymond,~\it La structure bor\'elienne d'Effros est-elle standard ?,~\rm 
Fund. Math.~100 (1978), 201-210

\end{document}